\documentclass[11pt]{amsart}
\usepackage{amsmath,amsfonts,amsthm,amssymb,
amsthm,graphicx,mathtools,amscd}
\usepackage{hyperref}

\usepackage[mathscr]{eucal}
\usepackage{graphics}
\usepackage{color}
\usepackage{epsfig}

\usepackage[nocompress]{cite}

\usepackage{rsfso}
\newtheorem{theorem}{Theorem}[section]
\newtheorem{lemma}[theorem]{Lemma}

\newtheorem{definition}[theorem]{Definition}

\newtheorem{remark}[theorem]{Remark}

\newcommand{\beginsec}{
\setcounter{equation}{0}
}

\newcommand{\la}{\lambda}
\newcommand{\eps}{\varepsilon}

\newcommand{\sig}{\sigma}
\newcommand{\del}{\delta}

\newcommand{\Del}{\mathnormal{\Delta}}

\newcommand{\La}{\mathnormal{\Lambda}}

\newcommand{\Sig}{\mathnormal{\Sigma}}

\newcommand{\N}{{\mathbb N}}

\newcommand{\R}{{\mathbb R}}

\newcommand{\Z}{{\mathbb Z}}

\newcommand{\E}{{\mathbb E}}

\newcommand{\PP}{{\mathbb P}}

\newcommand{\bS}{{\mathbf S}}

\newcommand{\scrA}{\mathscr{A}}
\newcommand{\scrC}{\mathscr{C}}
\newcommand{\scrD}{\mathscr{D}}
\newcommand{\scrG}{\mathscr{G}}

\renewcommand{\proof}{\noindent{\bf Proof.\ }}

\newcommand{\lan}{\langle}
\newcommand{\ran}{\rangle}

\newcommand{\w}{\wedge}

\newcommand{\To}{\Rightarrow}

\newcommand{\iy}{\infty}

\newcommand{\up}{\uparrow}

\newcommand{\cadlag}{c\`adl\`ag }
\newcommand{\noi}{\noindent}

\newcommand{\br}[1]{\langle #1 \rangle}
\newcommand{\sr}[1]{{\mathcal #1}}

\newcommand{\dd}[1]{\mathbb{#1}}
\newcommand{\eq}[1]{(\ref{eq:#1})}

\newcommand{\lem}[1]{Lemma~\ref{lem:#1}}

\newcommand{\app}[1]{Appendix~\ref{app:#1}}
\newcommand{\sectn}[1]{\S\ref{sec:#1}}

\newcommand{\sect}[1]{\ref{sec:#1}}

\usepackage{tcolorbox}
\newtcolorbox{bx}{colback=white!5!white,colframe=blue!75!black}

\begin{document}

\title[]{Heavy traffic limit with discontinuous coefficients via a non-standard semimartingale decomposition}

\author{Rami Atar}
\address{Viterbi Faculty of Electrical and Computer Engineering,
Technion -- Israel Institute of Technology, Haifa, Israel
}
\email{rami@technion.ac.il}
\author{Masakiyo Miyazawa}
\address{Department of Information Sciences, Tokyo University of Science, Noda, Chiba, Japan
}
\email{miyazawa@rs.tus.ac.jp}



\date{\today}

\begin{abstract}
This paper studies a single server queue in heavy traffic,
with general inter-arrival and service time distributions, where
arrival and service rates vary discontinuously as a function
of the (diffusively scaled) queue length. It is proved that
the weak limit is given by the unique-in-law solution to a stochastic
differential equation in $[0,\iy)$ with discontinuous drift
and diffusion coefficients.
The main tool is a semimartingale decomposition for point
processes introduced in \cite{dal-miy}, which is distinct from the 
Doob-Meyer decomposition of a counting process. Whereas the use of this tool
is demonstrated here for a particular model,
we believe it may be useful for investigating the scaling limits
of queueing models very broadly.
\end{abstract}

%\begin{quotation}
%\noindent {\bf keywords:} semimartingale decomposition; single server queue; level dependent arrival and service rates; heavy traffic; diffusion approximation
%
%\noindent {\bf AMS subject classification:} 60K25; 60H10; 60H30
%\end{quotation}

\maketitle

\section{Introduction}
\label{sec:introduction}

The goal of this paper is twofold: To continue a line of research
on the multi-level $GI/G/1$ queue at the diffusion scale
\cite{Miya2025, Miya2024b, Miya2024d},
%and to demonstrate the applicability of a semimartingale decomposition for general point processes, introduced in \cite{dal-miy}, to the study of queueing models under scaling limits.
and to demonstrate that a semimartingale decomposition for general point processes, introduced in \cite{dal-miy}, may serve as a powerful tool in the study of non-Markovian queueing models under scaling limits.

We consider a version of the multi-level $GI/G/1$ queueing model, which will be referred to in what follows as the {\it multi-level queue}.
In this model, arriving customers join an unlimited capacity buffer and are served in the first-come and first-served manner by a single server, and arrivals and services, driven by two independent renewal processes, undergo time change depending on the level to which the queue length belongs. Here, levels are disjoint intervals that form a finite partition of the half-line $\R_+=[0,\iy)$. In other words, the arrival and service rates
are determined by a piecewise constant function of the queue length.
Diffusion limits of models where arrival and service rates
vary with the queue length have been studied before
under continuous dependence; see \cite{wee14} and references therein.
However, there are very few limit results for models where the dependence is
discontinuous, because such discontinuities present some technical challenges.
Our main result is the weak convergence of the queue length,
normalized at the diffusion scale, to a reflected diffusion with
discontinuous drift and diffusion coefficients.

The proof of the result is based on a
semi-martingale decomposition for general point
processes derived in \cite{dal-miy}, referred to
as {\it the Daley-Miyazawa decomposition}.
This tool has recently been used in a crucial way
in \cite{ata-wol} to study a different non-Markovian queueing model in a different scaling regime, specifically, a load balancing model under
a mean-field many-server scaling. We believe that this decomposition,
along with some of its consequences developed in this paper in 
\S\ref{sec3}, may be useful in far greater generality
when it comes to proving scaling limits of queueing models.
From this viewpoint, this papers serves to provide an example
illustrating the applicability of the approach.

The Daley-Miyazawa decomposition for a counting process
was introduced in \cite{dal-miy} in
order to study extensions of Blackwell's renewal theorem.
%Note that a counting process is a special case of a nondecreasing process. Then,
This decomposition is distinct from the well-known
Doob-Meyer decomposition of counting processes (e.g., see Theorem 3.17 in Chapter I of \cite{jac-shi}), a special case of the Doob-Meyer decomposition
of submartingales (note that, as a nondecreasing process, a 
counting process is always a submartingale).
In the latter, the finite variation component is predictable, and called the {\it predictable compensator}. In the special case where this component is
absolutely continuous, its density is called
the {\it stochastic intensity} of the point process, a notion that
has been used very broadly in the study of point processes.
However, when it comes to queueing models, the stochastic intensity and,
more generally, the Doob-Meyer decomposition of counting processes, have rarely been
used beyond cases where the driving processes are
Poisson or compound Poisson. It seems that the reason for this is that
the stochastic intensity is not analytically tractable
beyond these cases. Contrary to this,
the finite variation component is not predictable in the Daley-Miyazawa decomposition. However, its first order term is equal
to the law-of-large-numbers limit of the point process,
a fact that makes it particularly attractive for scaling limit analysis.
Under diffusion scaling, the convergence of its martingale component
to a stochastic integral is based on the general theory of martingales,
and is quite different from the traditional approach to heavy traffic limits
which is based on the central limit theorem.

\subsection{Related work}

{\it On martingales associated with renewal processes.}
Martingales associated with renewal processes have been constructed 
before and used for analyzing
non-Markovian queues in heavy traffic \cite{kri89, kri-tak}.
To contrast these contributions with \cite{dal-miy}, we
point out several differences.
First, the Daley-Miyazawa decomposition addresses
general point processes, not only renewals. This is important
in the setting of this paper as seen in \S\ref{sec3} where we apply the 
decomposition to
complicated point processes denoted there by $A^n$ and $D^n$,
representing arrival and departure counts.
Further, it is not only the martingale component, but
the entire structure of the decomposition, including
the bounded variation component and the filtration, that comes into play
in the analysis. We have already mentioned a useful property
of the bounded variation component. As for the filtration,
it turns out that one can construct the decomposition of the various
processes on a common filtration (specifically, this is true
for $A^n$ and $D^n$, as we will show in Lemma \ref{lem1}),
a fact crucially used in our treatment. These aspects are intrinsic
to the decomposition proposed in \cite{dal-miy}
and are important elements of the approach.

{\it On convergence to diffusion with discontinuous coefficients.}
The paper \cite{kry-lip} proposed an approach to establishing
weak limits given by diffusion with discontinuous coefficients.
This approach gives general conditions in any dimension, based on
predictable characteristics. Although potentially applicable here,
the method seems less convenient to apply to queueing
systems driven by renewals. Indeed, the queueing model
provided in \cite{kry-lip} to demonstrate their approach is
a Markovian one.
The reader is referred to \cite{kry-lip} for some other works
proving convergence of queueing models
to diffusions with discontinuous coefficients.

One setting where such discontinuities arise naturally is that of
control problem formulations, where a cost is to be minimized
in the heavy traffic limit.
This was the case with \cite{ata-lev}, where a controlled queueing
problem gave rise to a diffusion with discontinuous drift,
and in \cite{ACR2,ACR1},
where, under various different settings, the limit process under
asymptotically optimal control was characterized as a diffusion
where both coefficients are discontinuous.

{\it On the multi-level queue.}
A model closely related to the one studied here
has appeared in \cite{Miya2025}, which is called the 2-level $GI/G/1$ queue.
The two models have the same level structure,
but the dynamics for arrivals is different. In the model from
\cite{Miya2025}, arrivals are subject to different renewal processes
according to the level to which the queue length belongs, while service speed is changed in the same way as in our model.
Our model simplifies the arrival structure by using a single renewal
process for arrivals. This simplification may be meaningful in application
because it may well reflect the influence due to congestion through
random time change. Theoretically, it enables us to derive a reflected
diffusion with discontinuous coefficients as a process limit, whereas only the weak convergence of the stationary distribution is studied in \cite{Miya2025}. Furthermore, as discussed there, its process limit may not be characterized by the form of (\ref{09}) except for a special case of its system parameters, so it would be much harder to derive the process limit for the model of \cite{Miya2025}. Thus, from theoretical viewpoint, the present result may be considered as the first step to solve this harder problem.

\subsection{Paper organization and notation}

This paper is structured as follows.
In the rest of this section, we introduce notation used throughout this paper. In \S\ref{sec:model-result}, the model is defined and the main
result is stated. A martingale toolbox based on the Daley-Miyazawa decomposition
is provided in \S\ref{sec3}. Finally, the proof of the main result
is given in \S\ref{sec:proof}.

Let $\N$ denote the set of natural numbers.
Let $\iota:[0,\iy)\to[0,\iy)$ denote the identity map. For $(\dd{S},d_\dd{S})$ a Polish space, let $\scrC(\R_+,\dd{S})$
and $\scrD(\R_+,\dd{S})$ denote the space of continuous and, respectively,
\cadlag paths $[0,\iy)\to\dd{S}$, endowed with the topology of uniform convergence
on compacts and, respectively, the $J_1$ topology (see \cite{Kall2001} for the details of this terminology). We typically put $\dd{S} = \dd{R}_{+} \equiv [0,\infty)$ or $\dd{S} = \dd{R}^{N}$, where $\R^N$ is the $N$-dimensional real vector space.
Denote by $\scrC^\up(\R_+,\dd{R}_{+})$ (resp., $\scrD^\up(\R_+,\dd{R}_{+})$)
the members of $\scrC(\R_+,\dd{R}_{+})$ (resp., $\scrD(\R_+,\dd{R}_{+})$)
that are nondecreasing and start at $0$.
Denote the Euclidean norm of $x \in \R^{N}$ by $\|x\|$.
For $\xi\in \scrD(\R_+,\R^N)$ and $T\in(0,\iy)$, denote
\begin{align*}
\Del\xi(t)&=\xi(t)-\xi(t-),\qquad t>0,\qquad \Del\xi(0)=\xi(0),
\\
w_T(\xi,\del)&=\sup\{\|\xi(t)-\xi(s)\|:s,t\in[0,T],|s-t|\le\del\},
\\
\|\xi\|^*_T&=\sup\{\|\xi(t)\|:t\in[0,T]\}.
\end{align*}
Denote by $W$ a standard Brownian motion in dimension $1$.
Denote by $\To$ convergence in distribution.
A sequence of probability measures on $\scrD(\R_+,\dd{S})$ is said to be $\scrC$-tight
if it is tight with the usual topology on $\scrD(\R_+,\dd{S})$ and every weak limit point
is supported on $\scrC(\R_+,\dd{S})$.
With a slight abuse of terminology, a sequence of random elements
(random variables or processes) is referred to as tight when their probability laws
form a tight sequence of probability measures, and a similar use is made for
the term $\scrC$-tight.

Some notation related to the Skorohod map on the half line is as follows
(see e.g.\ \cite[Section 8]{chu-wil} for details).
Given $\psi\in \scrD(\R_+,\R)$ with $\psi(0)\ge0$, there exists a unique
pair $(\phi,\eta)\in \scrD(\R_+,\dd{R}_{+})\times \scrD^\up(\R_+,\dd{R}_{+})$ such that
\begin{align*}
   \phi=\psi+\eta, \qquad \int_{[0,\iy)}\phi(t)d\eta(t)=0,
\end{align*}
and this pair is given by
\begin{equation}\label{19}
\eta(t)=\sup_{s\in[0,t]}\psi(s)^-,
\qquad
\phi(t)=\psi(t)+\eta(t),
\qquad
t\ge0,
\end{equation}
where $x^{-} = \max(0,-x)$ for $x \in \R$. The solution map from $\psi$ to $(\phi,\eta)$ is denoted by $\Gamma$, i.e., $\Gamma(\psi)=(\phi,\eta)$.

\section{Model and main results}
\label{sec:model-result}
\beginsec

In this section we
introduce a sequence of multi-level queues, and state heavy traffic conditions for its diffusion approximation. Then, we present a main theorem, whose proof will be given in the subsequent sections.

\subsection{Sequence of queueing models}
\label{sec:sequence}

We index the sequence of multi-level
queues by $n \in \N$. The number of threshold levels, denoted by $K \ge 2$, does not depend on $n$. For the $n$-th system, threshold levels are nonnegative numbers
denoted by $0=\ell^{n}_0<\ell^{n}_1<\cdots<\ell^{n}_{K-1}$, which partition $[0,\iy)$ into $K+1$ disjoint sets given as
\[
\bS^{n}_0=\{0\},\qquad \bS^{n}_i=(\ell^{n}_{i-1},\ell^{n}_i],\qquad 1\le i\le K-1,
\qquad \bS^{n}_K=(\ell^{n}_{K-1},\iy).
\]
We are given sequences of positive numbers
\begin{align}
\label{eq:la-mu}
  \la^{n}_{0}, \qquad \la^n_i, \quad \mu^n_i, \quad 1\le i \le K, \qquad n \in \N,
\end{align}
and $\mu^{n}_{0} = 0$. These numbers specify the inter-arrival rate and the service rate, respectively, at times when the queue length belongs to $\bS^{n}_{i}$. 

Service is assumed to be head-of-the-line and non-idling.
The queue length of the $n$-th system at time $t$, defined as the number of jobs in the system in that time, including the one being served (if there is any), is denoted by $X^{n}(t)$, and $X^{n} \equiv \{X^{n}(t); t \ge 0\}$ is referred to as the queue length process. Let
\begin{equation}\label{50}
I^n(t)=\int_0^t1_{\{X^n(s)=0\}}ds,
\qquad t\ge0,
\end{equation}
which is the cumulative idleness process. 
To define the arrival and service processes, we are given the two  mutually independent sequences of positive $i.i.d.$ random variables,
\begin{align}
\label{eq:ZAS}
  Z_A(j),\qquad Z_S(j),\qquad j=0,1,\ldots
\end{align}
having unit means and finite and positive variances $\sigma_{A}^{2}$ and $\sigma_{S}^{2}$, respectively. Then, define renewal processes $A$ and $S$ by
\begin{align}
\label{eq:AS-t}
  A(t) = \inf \Big\{i \ge 0; t < \sum_{j=0}^{i} Z_{A}(j)\Big\}, \quad S(t) = \inf \Big\{i \ge 0; t < \sum_{j=0}^{i} Z_{S}(j)\Big\}, \quad t \ge 0.
\end{align}
These notations without the index $n$ are commonly used for all the indexed systems.

We return to the $n$-th system. Let
\begin{align}
\label{eq:H-n}
  H^n_i(t)=\int_0^t1_{\{X^n(s)\in \bS^{n}_i\}}ds, \qquad 0\le i\le K,\ t\ge0,
\end{align}
and note that $H^n_0(t)=I^n(t)$. Further, let
\begin{align}
\label{eq:UVn-t}
  U^n(t)=\sum_{i=0}^K\la^n_iH^n_i(t),
\qquad
V^n(t)=\sum_{i=1}^K\mu^n_iH^n_i(t), \qquad t \ge 0.
\end{align}
Using $U^n(t)$ and $V^n(t)$, define the arrival and, respectively, departure counting processes of the $n$-th system by
\begin{equation}
\label{eq:AD-n}
A^n(t)=A(U^n(t)),
\qquad
D^n(t)=S(V^n(t)), \qquad t \ge 0.
\end{equation}
Then we have
\begin{equation}
\label{eq:X-n}
X^n(t) = X^{n}(0) + A^n(t)-D^n(t), \qquad t \ge 0.
\end{equation}
For simplicity it will be assumed that $X^n(0)=0$.
We regard the processes $A$ and $S$ as the stochastic primitives,
and the tuple $(X^n,I^n,(H^n_i)_{0\le i\le K},U^n,V^n,A^n,D^n)$,
which represents the queueing dynamics, to be determined
from the primitives via the system of equations
\eqref{50}, \eqref{eq:H-n}--\eqref{eq:X-n}.
This system of equations
has a unique solution given the primitives, as can
be proved by induction on the jump times. The inductive argument
is standard (e.g., see \cite[Section 2.5]{ChenMand1991}), and we omit the details.
Note that, by construction, $A$ and $S$ are right-continuous,
and as a result so are $X^n,A^n,D^n$.

\subsection{Heavy traffic conditions and main theorem}
\label{sec:heavy traffic}

We now introduce assumptions on the scaling of the various parameters.
The first assumption is about the levels. Assume that, for given constants  $\ell_{0} = 0 < \ell_{1} < \ldots < \ell_{K-1}$, the levels are given by
\begin{align}
\label{eq:ell-n}
\ell^{n}_0=0,
\qquad
  \ell^{n}_{i} = n^{1/2} (\ell_{i}+o(1)),
  \qquad
  1 \le i \le K-1,
\end{align}
where we write $f(n) = o(1)$ for function $f: \dd{N} \to \dd{R}$ if $f(n)$ vanishes as $n \to \infty$. Using these $\ell_{i}$, define a partition of $\dd{R}_{+}$ as
\begin{align*}
  \bS_0=\{0\},\qquad \bS_i=(\ell_{i-1},\ell_i],\qquad 1\le i\le K-1,
\qquad \bS_K=(\ell_{K-1},\iy).
\end{align*}
The second assumption is about the state dependent arrival and service rates given by \eq{la-mu}. Let constants
\[
\la_0\in(0,\iy), \qquad
\la_{i}, \mu_i\in(0,\iy), \quad \hat{\la}_{i}, \hat{\mu}_{i}\in\R, \quad 1\le i\le K,
\]
be given. Recall that $\la^{n}_{i}$ for $1\le i\le K$ and $\mu^{n}_i\in(0,\iy)$ for $1\le i\le K$ are given for the $n$-th system. Let
\begin{align}
\label{eq:bar la-mu}
 & \bar\la^n_i:=n^{-1}\la^n_i,
\qquad
0\le i\le K,
\qquad
\bar\mu^n_i:=n^{-1}\mu^n_i.
\qquad
1\le i\le K,\\
\label{eq:hat la-mu}
 & \hat\la^n_i:=n^{-1/2}(\la^n_i-n\la_i),
\qquad
\hat\mu^n_i:=n^{-1/2}(\mu^n_i-n\mu_i)
\qquad
1\le i\le K.
\end{align}
Assume that 
\begin{align}
\label{eq:la-0}
  \bar \la^{n}_{0} \to \la_{0}, \qquad n \to \infty,
\end{align}
and that, as $n \to \infty$,
\begin{align}
\label{eq:la-n}
 & \hat \lambda^{n}_{i}
 \to \hat \lambda_{i}, \qquad
 \hat \mu^{n}_{i}
 \to \hat{\mu}_{i}, \qquad 1\le i\le K.
\end{align}
Note that, from \eq{bar la-mu} and \eq{hat la-mu}, we have
\begin{align}\label{c01}
  \la^{n}_{i} = n \la_{i} + n^{1/2} (\hat \la_{i} + o(1)), \qquad  \mu^{n}_{i} = n \mu_{i} + n^{1/2} (\hat \mu_{i} + o(1)), \qquad 1\le i\le K.
\end{align}
Let the critical load condition hold, namely
\[
\la_i=\mu_i,\qquad 1\le i\le K.
\]

Let
\begin{align}
\label{hat-XI}
\hat{X}^{n}(t) = n^{-1/2} X^{n}(t),
\qquad
\hat I^n(t)=n^{1/2}\bar\la^n_0I^n(t),
\qquad t\ge0,
\end{align}
and let $b^n_i=\hat\la^n_i-\hat\mu^n_i$, $b_i=\hat\la_i-\hat\mu_i$, $\sig_i=(\la_i\sig_A^2+\mu_i\sig_S^2)^{1/2}$ for $1\le i\le K$, and $\sig_0=\la_0^{1/2}\sig_A$. Define
\begin{equation}\label{11}
b(x)=\sum_{i=1}^Kb_i1_{\bS_i}(x),
\qquad
\sig(x)=\sum_{i=0}^K\sig_i1_{\bS_i}(x), \qquad x \in \dd{R}_{+}.
\end{equation}
Consider the SDE for $(X,L)\in
\scrC(\R_+,\R_+)\times \scrC^\up(\R_+,\dd{R}_{+})$
\begin{equation}\label{09}
\begin{split}
&X(t)=x_0+\int_0^tb(X(s))ds+\int_0^t\sig(X(s))dW(s)+L(t) , \qquad t\ge0,
\\
&\int_{[0,\iy)}X(t)dL(t)=0.
\end{split}
\end{equation}
All assumptions made thus far are in force throughout the paper.
The main result of this paper is the following.

\begin{theorem}\label{th1}
There exists a weak solution to \eqref{09}, and it is
unique in law. Moreover, denoting by $(X,L,W)$ a solution 
corresponding to $x_0=0$,
and by $(\hat X^n,\hat I^n)$ the processes associated with
the mutli-level queue as in \eqref{hat-XI},
one has $(\hat X^n,\hat I^n)\To(X,L)$
in $\scrD(\R_+,\dd{R}_{+}) \times \scrC^\up(\R_+,\dd{R}_{+})$.
\end{theorem}

Theorem \ref{th1} will be proved in \sectn{proof}. For this proof, we prepare our key tool in the next section.

\begin{remark}[About our assumptions]
Inter-arrival and service times are not assumed to be continuously distributed, hence arrival and service completion may occur simultaneously. This presents no difficulty to the approach, but requires a computation of the quadratic cross-variation between martingale components of the arrival and departure processes (see Lemma~\ref{lem1}).

The assumption that the system starts empty is posed only to simplify the construction of the model. Under a more general initial system state, the limit process would involve a non-zero initial condition in \eqref{09}.
\end{remark}

\begin{remark}[About extensions]
It is well known that strong existence for SDE of the type \eqref{09} fails when $\sigma$ is discontinuous (see e.g.\  \cite[p.\ 375]{rev-yor}).
Closely related is the fact that discrete approximations of the dynamics do not, in general, converge to solutions of the SDE in the presence of such discontinuities. As pointed out by one of the referees, one may ask whether the main result can be extended to the case where the sequences of inter-arrival and potential service processes only satisfy a functional central limit theorem. Whereas in the case of Lipschitz coefficients it is well known that the answer is positive and follows from the continuous mapping theorem, the answer here is very likely negative, for the reason mentioned above.

An extension of our setting that is left for future work is to network models, with level-dependent service and arrival intensities, which may lead to SDE in higher dimension with discontinuities along surfaces.

Another extension yet to be explored is the multi-level queueing model mentioned in the introduction. When the arrivals are driven by different renewal processes at different levels, there is a subtlety related to the behavior of the queue length and the time spent near the points of discontinuity, leading us to believe that the limiting SDE would be a variant of \eqref{09} that includes additional local time terms for the time spent at points of discontinuity. It is likely that the Daley--Miyazawa tool will be useful in this setting as well.
\end{remark}

\section{Daley-Miyazawa decompositions}
\label{sec3}
\beginsec

Recall that $A$ and $S$ were constructed based on
$Z_A(j), Z_S(j)$, $j=0,1,\ldots$,
denoting the times between successive counting events. Let
\[
\zeta_A(j)=1-Z_A(j),\qquad \zeta_S(j)=1-Z_S(j), \qquad j=0,1,\ldots
\]
Let
\[
R_A(t)=\inf\{s>0:A(t+s)>A(t)\}, \qquad t\ge0
\]
denote the residual time to the next counting event. Define similarly
$R_S(t)$. By definition,
\begin{equation}\label{05}
t+R_A(t)=\sum_{j=0}^{A(t)}Z_A(j),
\qquad
t+R_S(t)=\sum_{j=0}^{S(t)}Z_S(j).
\end{equation}
Denote
\begin{equation}\label{06}
M_A(t)=\sum_{j=1}^{A(t)}\zeta_A(j),
\qquad
M_S(t)=\sum_{j=1}^{S(t)}\zeta_S(j).
\end{equation}
We then have the Daley-Miyazawa semimartingale representation
\cite{dal-miy}
\begin{equation}\label{07}
A(t)=t+R_A(t)-Z_A(0)+M_A(t),
\qquad
S(t)=t+R_S(t)-Z_S(0)+M_S(t),
\qquad t\ge0,
\end{equation}
with $M_A$ and $M_S$ being martingales with respect to the filtrations generated, respectively, by $\{(A(t),R_{A}(t)); t \ge 0\}$ and $\{(S(t),R_{S}(t)); t \ge 0\}$. Our main interest will be in
martingales that can be viewed as time changed versions of the
above ones.
These are
\begin{equation}\label{08}
\hat M^n_A(t)=n^{-1/2}\sum_{j=1}^{A^n(t)}\zeta_A(j),
\qquad
\hat M^n_S(t)=n^{-1/2}\sum_{j=1}^{D^n(t)}\zeta_S(j),
\qquad
\hat M^n=\hat M^n_A-\hat M^n_S.
\end{equation}
Note, by \eqref{07}, that, for $t\ge0$,
\begin{equation}\label{34}
\begin{split}
A^n(t)=A(U^n(t))=U^n(t)+R_A(U^n(t))-Z_A(0)+n^{1/2}\hat M^n_A(t),
\\
D^n(t)=S(V^n(t))=V^n(t)+R_S(V^n(t))-Z_S(0)+n^{1/2}\hat M^n_S(t).
\end{split}
\end{equation}
Of crucial importance to the analysis will be the fact, proved next,
that \eqref{34} gives a semimartingale decomposition
of the arrival and departure processes, $A^n$ and $D^n$, on a {\it single} filtration.
To make this precise, let the `state' of the system be given by
\[
\sr{S}^n(t)=(X^n(t),A^n(t),D^n(t),R_A(U^n(t)), R_S(V^n(t)),\tilde R^n_A(t),\tilde R^n_D(t)),
\qquad t\ge0,
\]
where
\[
\tilde R^n_A(t)=\inf\{s>0:A^n(t+s)>A^n(t)\},
\qquad
\tilde R^n_D(t)=\inf\{s>0:D^n(t+s)>D^n(t)\},
\qquad
t\ge0
\]
denote the residual times to the next counting events
of $A^n$ and $D^n$.
Let
\[
\sr{F}^n_t=\sig\{\sr{S}^n(s):0\le s\le t\},
\qquad t\ge0.
\]
Because of the inclusion of $\tilde R^n_A$ and $\tilde R^n_D$,
the (right-continuous) counting processes $A^n$
and $D^n$ are $\{\sr{F}^n_t\}$-predictable.

We next compute the (predictable and optional) quadratic (cross-) variations associated with $\hat M^n_A$ and $\hat M^n_S$.
%their cross variations and their optional ones, where the cross variation $[\hat M^n_A,\hat M^n_S]$
Recall that if $\xi_1$ and $\xi_2$ are piece-wise constant semimartingales then their optional cross-variation can be computed as
\begin{equation}
\label{g10}
[\xi_1,\xi_2](t)=\sum_{0<s\le t}\Delta \xi_1(s)\Delta\xi_2(s),\qquad t\ge0,
\end{equation}
whereas the optional quadratic variation is given by $[\xi_1](t)=[\xi_1,\xi_1](t)$
%\[
%  [\hat M^n_A, \hat M^n_S](t) = \Big(\sum_{0<s\le t} \Delta \hat M^n_A(s) \hat M^n_S(s)\Big), \qquad t \ge 0,
%\]
%because $\hat M^n_A(t)$ and $\hat M^n_S(t)$ are piecewise constant functions 
(see \cite[Definition I.4.45, Theorem I.4.47]{jac-shi} for the definition of cross-variation and its computation).
The predictable cross-variation (resp., predictable quadratic variation), denoted by $\br{\xi_1,\xi_2}$ (resp., $\br{\xi_1}$) is defined as the predictable component of $[\xi_1,\xi_2]$ (resp., $[\xi_1]$).
% Similar to an optional quadratic variation, the optional quadratic variation $\br{\hat M^n_A,\hat M^n_S}$ is defined as the predictable component of $[\hat M^n_A,\hat M^n_S]$.

\begin{lemma}\label{lem1}
The processes
$\hat M^n_A$ and $\hat M^n_S$ are $\{\sr{F}^n_t\}$-square
integrable martingales, with optional quadratic variations
\begin{equation}\label{40}
[\hat M^n_A](t)=n^{-1}\sum_{j=1}^{A^n(t)}\zeta_A(j)^2
\qquad
[\hat M^n_S](t)=n^{-1}\sum_{j=1}^{D^n(t)}\zeta_S(j)^2,
\qquad t\ge0,
\end{equation}
and the optional cross-variation $[\hat M^n_A,\hat M^n_S](t)$,
is an $\{\sr{F}^n_t\}$-martingale.
Moreover, the predictable quadratic variations are given by
\begin{equation}\label{42}
\lan\hat M^n_A\ran(t)=n^{-1}\sig^2_AA^n(t)
\qquad
\lan\hat M^n_S\ran(t)=n^{-1}\sig^2_SD^n(t),
\qquad
\lan\hat M^n_A,\hat M^n_S\ran(t)=0,
\qquad
t\ge0.
\end{equation}
\end{lemma}

\begin{remark}
If either interarrival times or service times have continuous distribution, arrivals and departures never occur at the same time, a.s. In this case the proof somewhat simplifies, as $[\hat M^n_A,\hat M^n_S](t)=0$ a.s. However, this is not true in general.
\end{remark}

\proof
Because $X^n$ is, by definition, $\{\sr{F}^n_t\}$-adapted,
so are $H^n_i$, $U^n$ and $V^n$. Because
$R_A(U^n(\cdot))$ is $\{\sr{F}^n_t\}$-adapted and
$R_A(U^n(0))=R_A(0)=Z_A(0)$, one has that $Z_A(0)$ is
$\sr{F}^n_0$-measurable. The same is true for $Z_S(0)$.
In view of the identities \eqref{34},
$\hat M^n_A$ and $\hat M^n_S$ are also $\{\sr{F}^n_t\}$-adapted.

We turn to proving square integrability.
To this end, note first that, as a renewal process, $\E A(t)<\iy$
for all $t$. Since, by \eqref{c01}, $U^n(t)\le c n t$ for some constant $c$, it follows that
$\E A^n(t)<\iy$ as well, and by a similar argument, $\E D^n(t)<\iy$
for all $t$. Next, for an $\{\sr{F}^n_t\}$-stopping time $\tau$, denote
\begin{equation}\label{35-}
\sr{F}^n_{\tau-}
=\sr{F}^n_0\vee\sig\{C\cap\{t<\tau\}: C\in \sr{F}^n_t, t\ge0\}.
\end{equation}
Then it is a standard fact that $\tau\in\sr{F}^n_{\tau-}$
\cite[I.1.11, I.1.14]{jac-shi}. Consider
\[
s^n(k)=\inf\{t\ge0: A^n(t)\ge k\},
\qquad
t^n(k)=\inf\{t\ge0: D^n(t)\ge k\},
\qquad k=1,2,\ldots
\]
These are $\{\sr{F}^n_t\}$-stopping times, and as a consequence
\begin{equation}\label{35}
s^n(k)\in\sr{G}^n_k:=\sr{F}^n_{s^n(k)-},
\qquad
t^n(k)\in\sr{H}^n_k:=\sr{F}^n_{t^n(k)-},
\qquad k=1,2,\ldots
\end{equation}
We now argue that
\begin{equation}\label{36}
\text{$\zeta_A(j)\in\sr{G}^n_k$  for $j\le k-1$, and
$\{\zeta_A(l),l\ge k\}$ is independent of $\sr{G}^n_k,\qquad k=1,2,\ldots$}
\end{equation}
and
\begin{equation}\label{37}
\text{for every $t$, $A^n(t)$ is a stopping time on the discrete parameter
filtration $\{\sr{G}^n_k\}$.}
\end{equation}
To this end, let $k\ge1$. The first assertion in \eqref{36} follows from the fact that
$\{Z_A(j), j\le k-1\}$ can be determined from $\{s^n(j), j\le k\}$
and $\{U^n(s), s<s^n(k)\}$.
We next consider the tuple
\[
\Theta_k=\{Z_A(j),\, j=0,1,\ldots, k-1, Z_{S}(j),\, j=0,1,2,\ldots\}.
\]
By the construction of the model, all the information about the state of the system
up to $s^n(k)-$, namely $\{\sr{S}^n(t),t<s^n(k)\}$, is determined by $\Theta_{k} \in \{\sr{G}^n_k\}$,
whereas $\{Z_A(l), l\ge k\}$ (hence $\{\zeta_A(l),l\ge k\}$)
is independent of $\Theta_k$. This shows the second assertion in \eqref{36}.
Thus, \eqref{36} is proved.
Next, by \eqref{35},
\[
\{A^n(t)\le k\}=\{s^n(k)\ge t\}\in\sr{G}^n_k,
\]
proving \eqref{37}.

The structure \eqref{35}--\eqref{36} makes it possible to apply Wald's equations.
In particular, by Wald's second equation \cite[Theorem 4.1.6]{dur-bk},
\[
\E\Big[\Big(\sum_{j=1}^{A^n(t)}\zeta_A(j)\Big)^2\Big]
=\E[A^n(t)]\E[\zeta_A(0)^2]<\iy.
\]
This shows that $\E[\hat M^n_A(t)^2]<\iy$ for all $t$.
The argument for $\hat M^n_S$ is completely analogous.

For the martingale property of $\hat M^n_A$, note first by \eqref{36} that
\begin{equation}\label{38}
\E[\zeta_A(k)|\sr{G}^n_k]=0.
\end{equation}
Next, as in the proof of \cite[Lemma 2.1]{dal-miy}, write
\[
\hat M^n_A(t)=n^{-1/2}\sum_{k=1}^{A^n(t)}\zeta_A(k)
=n^{-1/2}\sum_{k\ge1}\zeta_A(k)1_{\{s^n(k)\le t\}},
\]
and for $s<t$,
\begin{align}\label{38+}
\E[\hat M^n_A(t)|\sr{F}^n_s]-\hat M^n(s)
&= n^{-1/2}\sum_{k\ge1}\E[\zeta_A(k)1_{\{s<s^n(k)\le t\}}|\sr{F}^n_s].
\end{align}
By \eqref{35},
\begin{equation}\label{ww1}
\{s<s^n(k)\le t\}\in \sr{G}^n_k.
\end{equation}
By \eqref{35-}, if $C\in\sr{F}^{n}_{s}$ then $C\cap\{s<s^n(k)\}\in\sr{G}^{n}_{k}$. As a result, by (\ref{38}),
\begin{align*}
  \dd{E}[\zeta_A(k) 1_{\{s<s^n(k)\le t\}} 1_{C}] & = \dd{E}[\dd{E}[\zeta_A(k) 1_{\{s<s^n(k)\le t\}} 1_{C}|\sr{G}^{n}_{k}]] \\
  & = \dd{E}[\dd{E}[\zeta_A(k)|\sr{G}^{n}_{k}] 1_{\{s<s^n(k)\le t\}} 1_{C}] = 0, \qquad \forall C \in \sr{F}^{n}_{s}.
\end{align*}
This gives
\begin{equation}\label{ww3}
\dd{E}[\zeta_A(k) 1_{\{s<s^n(k)\le t\}}|\sr{F}^{n}_{s}] = 0.
\end{equation}
Hence,
going back to \eqref{38+}, we have
 \begin{align*}
\E[\hat M^n_A(t)|\sr{F}^n_s]-\hat M^n(s)
&=0.
\end{align*}
A similar argument holds for $\hat M^n_S$.
This completes the proof that $\hat M^n_A$ and $\hat M^n_S$
are square integrable $\{\sr{F}^n_t\}$-martingales.

The expression stated for $[\hat M^n_A](t)$ in \eqref{40}
(similarly for $[\hat M^n_S](t)$) is immediate from
(\ref{08}), which can be written as $\hat M^n_A(t) = \sum_{s\le t}\Del\hat M^n_A(s)$ because it is piecewise
constant with finitely many jumps on finite intervals.

As for the cross-variation, by \eqref{g10},
\begin{align*}
m^n(t):=[\hat M^n_A,\hat M^n_S](t)
&=n^{-1}\sum_{k=1}^{A^n(t)}\sum_{l=1}^{D^n(t)}\zeta_A(k)\zeta_S(l)1_{\{s^n(k)=t^n(l)\}}, \qquad t \ge 0.
\end{align*}
The $\{\sr{F}^n_t\}$-adaptedness of $m^n$ is clear from that of
$\hat M^n_A$ and $\hat M^n_S$, and absolute integrability of the former
follows from square integrability of the latter.
Let now $u^n(k,l)=\min(s^n(k),t^n(l))$. For $k,l\ge1$, let
\[
\Theta_{k,l}=\{Z_A(k'),\,k'<k,\,, Z_S(l'),\,l'<l\}.
\]
Then, by construction, all information about $\{\sr{S}(s), s<u^n(k,l)\}$ can be
determined by $\Theta_{k,l}$, whereas $Z_A(k)$, $Z_S(l)$ and $\Theta_{k,l}$
are mutually independent. This gives
\begin{equation}\label{44}
\E[\zeta_A(k)\zeta_S(l)|\sr{F}^n_{u^n(k,l)-}]=0,
\qquad
k,l\ge1.
\end{equation}
%\\
%&=
%n^{-1}\sum_{k\ge1}\sum_{l\ge1}
%\E[\E[\zeta_A(k)\zeta_S(l)1_{\{s<s^n(k)=t^n(l)\le t\}}| \sr{F}^n_{u^n(k,l)-}]|\sr{F}^n_s],
Next, for $0\le s<t$, it follows from \eqref{35} that
\begin{equation}\label{ww4}
\{s<s^n(k)=t^n(l)\le t\}\in\sr{F}^n_{u^n(k,l)-}.
\end{equation}
Using the same reasoning by which \eqref{ww3} was derived from \eqref{38} and \eqref{ww1}, we can now deduce
\[
\E[\zeta_A(k)\zeta_S(l)1_{\{s<s^n(k)=t^n(l)\le t\}}| \sr{F}^n_s]=0
\]
from \eqref{44} and, respectively, \eqref{ww4}. Thus
\begin{align*}
\E[m^n(t)|\sr{F}^n_s]-m^n(s)&=n^{-1}\sum_{k\ge1}\sum_{l\ge1}
\E[\zeta_A(k)\zeta_S(l)1_{\{s<s^n(k)=t^n(l)\le t\}} |\sr{F}^n_s]=0,
\end{align*}
completing the proof of the martingale property of $m^n$.

For the first part of \eqref{42}, it suffices to show that
$n^{-1}\sig^2_AA^n(t)$ is the compensator of $[\hat M^n_A](t)$
(see \cite[Proposition I.4.50 p.\ 53]{jac-shi}). This is an immediate consequence of $A^n$ being $\{\sr{F}^n_t\}$-predictable
and $\E[\zeta_A(j)^2] =\sig^2_A$ (indeed, the jumps are predictable because the residual times are included in the filtration). A similar proof holds for
the expression $\lan\hat M^n_S\ran$.
Finally, $\lan\hat M^n_A,\hat M^n_S\ran$ is the compensator of
$[\hat M^n_A,\hat M^n_S]$. Since the latter is a zero mean martingale,
it follows that the former is zero.
\qed

\section{Proof of main result}
\label{sec:proof}
\beginsec

The proof of Theorem \ref{th1} is established via a compactness--uniqueness argument. In particular, it is shown that uniqueness in law holds for equation \eqref{09} (Lemma \ref{lem01}), that the rescaled processes $(\hat X^n,\hat I^n)$ form a tight sequence (Lemma \ref{lem02}), and that every subsequential limit  of this pair forms a solution to \eqref{09} (Lemma \ref{lem03}). To state these lemmas, some additional terminology and notation are required. We prepare them in the next subsection.

\subsection{Preliminaries}
\label{sec:preliminary}

We first recall some fundamental facts about semimartingales and the solution of the SDE \eqref{09}.

\begin{itemize}
\item [(i)] For a semimartingale $\xi$, $\Del [\xi] = (\Del \xi)^{2}$, and therefore, the sample paths of $[\xi]$ are continuous a.s.\ if and only if $\xi$ is a continuous semimartingale.
\item [(ii)] For a continuous semimartingale $X$ and each $a \in \dd{R}$, one of the equivalent definitions of the local time $L^{X}_{a}(t)$ of $X$ at $a$ is
\begin{align}
\label{eq:local-t}
  L^{X}_{a}(t) = \lim_{\varepsilon \downarrow 0} \frac 1{\varepsilon} \int_{0}^{t} 1_{\{X(s) \in [a,a+\eps)\}} d\br{X}, \qquad t \ge 0.
\end{align}
Moreover, $L^{X}_{a}$ exists $a.s.$, and has a version that it is continuous in $t$ and continuous in $a$ from the right with left-hand limits. 
\item [(iii)] For any solution $(X,L)$ to the SDE \eqref{09}, the boundary term $L$ and the local time $L^X_0(t)$ of $X$ at $0$ are related via $L(t)=\frac{1}{2}L^X_0(t)$ a.s.
\end{itemize}

The proofs of these three facts can be found in textbooks. For example, (i) is \cite[Theorem I.4.47 (c)]{jac-shi}. For (ii) and (iii), see \cite[Corollary VI.1.9]{rev-yor} and \cite[Section 1.3]{pil14}, respectively. By these facts, the SDE \eqref{09} is equivalent to the following SDE for $X \in \scrC(\R_+,\dd{R}_{+})$.
\begin{equation}\label{eq:ref-d}
  X(t)=x_0+\int_0^tb(X(s))ds+\int_0^t\sig(X(s))dW(s)+\frac 12 L^{X}_{0}(t) , \qquad t\ge0,
\end{equation}

We aim at showing that there is a unique weak solution $(X,W)$ of the SDE \eq{ref-d}, equivalently, \eqref{09}.
% However, this is not obvious because $a(x)$ and $\sigma(x)$ are discontinuously changed.
Although such results for SDE on $\R$ (and on $\R^N$)
with discontinuous coefficients are well-known, there are very few results concerning SDE with reflection and discontinuous coefficients. To make our exposition self-contained, we include a complete proof.
To this end, we consider another way to characterize the solution $X$, referred to as the {\it martingale problem}, known as a powerful tool for constructing Markov processes \cite{eth-kur,Kall2001,StroVara1979}. It has different versions depending on whether $X$ has a boundary. We introduce two of them for a continuous semimartingale $X$. For this, we use the following notation.

Let ${\scrC}^{\infty}_0(\dd{S},\dd{R})$ be the set of compactly supported functions $\dd{S} \to \R$ that are infinitely differentiable, where $\dd{S} = \dd{R}$ or $\dd{S} = \dd{R}_{+}$, and the differentiability at the origin is from right for $\dd{S} = \dd{R}_{+}$.
% Note that $\scrC^{\infty}_0(\dd{R}_{+},\dd{R})$ is a subset of $\scrC(\R_+,\R)$, but they are used for different objects, a set of test functions and the set of continuous sample paths, respectively.

We first recall the martingale problem for a semimartingale $X$ with sample paths in $\scrC(\R_+,\R)$ corresponding to an SDE on all of $\R$. While the notion is standard, we mention it here in order to build on it the corresponding notion for SDE with reflection, such as \eqref{eq:ref-d}. Let $\scrG$ be an operator on ${\scrC}^{\infty}_{0}(\R,\R)$. For each probability distribution $\nu$ on $\dd{R}$, if $x_0 \sim \nu$ and
\begin{align}
\label{eq:g-mp}
  f(X(t))- \int_0^t \scrG f(X(s)) ds \mbox{ is a local martingale for all } f \in {\scrC}^{\infty}_{0}(\R,\dd{R}),
\end{align}
then $X$ is called a solution of the martingale problem for $(\nu,\scrG,\dd{R})$ (see \cite[Section 2.3]{eth-kur} for the definition of a local martingale). When $\scrG$ is the elliptic operator, that is,
\begin{align}
\label{eq:G}
  \scrG f (x) = u(x) f'(x) + \frac 12 v(x) f''(x), \qquad f \in {\scrC}^{\infty}_{0}(\R,\dd{R}),
\end{align}
for measurable bounded functions $u, v$ from $\dd{R}$ to $\dd{R}$ such that $v(x) > 0$ for $x \in \dd{R}$, $X$ is the diffusion process with drift $u$ and diffusion coefficient $v$ subject to the initial distribution $\nu$ if and only if $X$ is the solution of the martingale problem $(\nu,\scrG,\dd{R})$. The proof of this important characterization can be found in text books (e.g., see \cite[Theorem 21.7]{Kall2001}).

We next modify this notion so as to characterize the reflecting process $X$ of our interest. There are several ways to do this. Among the ones that are mentioned in \cite[Section 2]{StroVara1971},
% although the diffusion coefficient is continuous there.
we take the martingale condition in \cite[Theorem 2.4]{StroVara1971}, and refine it by (iii). Consider the operator on the domain ${\scrC}^{\infty}_{0}(\R_+,\R)$ defined as
\begin{align}
\label{eq:A}
  \scrA f(x)=b(x)f'(x)+\frac{1}{2}\sig^2(x)f''(x), \qquad x \in \R, f \in {\scrC}^{\infty}_{0}(\R_+,\R),
\end{align}
with $b,\sig$ as in \eqref{11}.

\begin{definition}
\label{dfn:e-mp}
A semi-martingale $X$ with sample paths in $\scrC(\R_+,\R_{+})$ is said to solve the martingale problem for $(\nu,\scrA,\R_{+})$ if $x_0 \sim\nu$ and, denoting $\dd{F}^{X} = \{\sr{F}^{X}_{t}; t \ge 0\}$,
\begin{align}
\label{eq:e-mp}
\begin{split}
  f(X(t))- & \Big(\int_0^t\scrA f(X(s)) ds + \frac 12 \int_{0}^{t} f'(X(s)) dL^{X}_{0}(s)\Big) \\
  & \qquad \mbox{ is an $\dd{F}^{X}$-local martingale for all } f \in {\scrC}^{\infty}_{0}(\R_+,\R).
\end{split}
\end{align}
\end{definition}

The next lemma shows that the solution of the SDE \eq{ref-d} is indeed characterized by this martingale problem.
\begin{lemma}
\label{lem:e-mp}
Let $X$ be a continuous semimartingale. Then $X$ is a weak solution of the SDE \eq{ref-d} if only if $X$ is a solution of the martingale problem for $(\nu,\scrA,\dd{R}_{+})$. Furthermore, both solutions are unique in law if either one of them is unique in law.
\end{lemma}

This lemma can be proved in the exactly same way as the proof of \cite[Theorem 21.7]{Kall2001}. However, we detail its proof in \app{e-mp} because of its importance in our proofs.

\subsection{Key lemmas and proof of Theorem \ref{th1}}
\label{key-lem}

We present key Lemmas \ref{lem01}--\ref{lem03}, from which Theorem \ref{th1} follows. We first consider the unique existence of a solution of the SDE \eqref{09}. 

\begin{lemma}\label{lem01}
Let $\nu$ be a Borel probability measure on $[0,\iy)$. Then there exists a unique solution to
the martingale problem for $(\nu,\scrA,\dd{R}_{+})$. Equivalently, there exists
a unique weak solution $(X,L,W)$ to SDE \eqref{09} (namely, $(X,W)$ to SDE \eq{ref-d}), with $x_0\sim\nu$,
and the law of $X$ is equal to the unique solution to the martingale problem $(\nu,\scrA,\dd{R}_{+})$.
\end{lemma}

This lemma will be proved in \sectn{lem01}. It is closely related to \cite[Lemma 2.2]{Miya2024b}, but Lemma 2.2 there is nothing to do with any martingale problem. It is also very close to \cite[Lemma 4.1(3)]{ACR2}.

Next, note that by \eqref{08} and \eqref{34}, we can write \eqref{eq:X-n} as
\[
\hat X^n=n^{-1/2}X^n=\hat M^n+n^{-1/2}(U^n-V^n)+e^n,
\]
where
\begin{equation}\label{13}
e^n(t)=n^{-1/2}\Big[R_A(U^n(t))
-R_S(V^n(t))-Z_A(0)+Z_S(0)\Big],
\qquad t\ge0.
\end{equation}
Moreover, by \eqref{eq:bar la-mu}, \eqref{eq:hat la-mu} and the definition of $b^n_i$, we have $n^{-1/2}(\la^n_i- \mu^n_i)=b^n_i$ for $1\le i\le K$, and $n^{-1/2} \lambda^{n}_{0} = n^{1/2} \bar{\lambda}^{n}_{0}$. Therefore
\[
n^{-1/2}(U^n-V^n)
=\sum_{i=1}^Kb^n_iH^n_i+\hat I^n,
\]
where we have used \eqref{eq:UVn-t}, the definition \eqref{hat-XI} of $\hat I^n$, recalling from \eqref{eq:H-n} that $H^n_0=I^n$.
We arrive at the relation
\begin{equation}\label{10}
\hat X^n=\hat Y^n+\hat I^n
\qquad
\text{where}
\qquad
\hat Y^n=\hat M^n+\sum_{i=1}^Kb^n_iH^n_i+e^n.
\end{equation}
Finally, let $\bar A^n=n^{-1}A^n$ and $\bar D^n=n^{-1}D^n$, and
\[
e^n_A=[\hat M^n_A]-\lan\hat M^n_A\ran,
\qquad
e^n_{S}=[\hat M^n_S]-\lan\hat M^n_S\ran.
\]

\begin{lemma}\label{lem02}
The tuple
\begin{equation}\label{16}
\Sig^n:=(\bar A^n,\bar D^n,\hat X^n,\hat I^n,\hat M^n,[\hat M^n],\hat Y^n, (H^n_i)_{0\le i\le K}, e^n, e^n_A, e^n_S)
\end{equation}
is $\scrC$-tight in $\scrD(\R_+,\R)^{10+K}$, and every subsequential weak limit is given by
\begin{equation}\label{18}
\Sig:=(\bar A,\bar D,X,L,M,[M],Y,(H_i)_{0\le i\le K},0,0,0),
\end{equation}
where $M$ is a martingale, and the following relations hold:
\begin{equation}\label{20-}
(X,L)=\Gamma(Y),
\qquad
H_0=0,
\qquad
Y=M+\sum_{i=1}^Kb_iH_i,
\end{equation}
\begin{equation}\label{20}
\bar A=\sum_{i=1}^K\la_iH_i,
\qquad
\bar D=\sum_{i=1}^K\mu_iH_i,
\qquad
[M]=\sum_{i=1}^K\sig^2_iH_i.
\end{equation}
\end{lemma}

\begin{lemma}\label{lem03}
Given a subsequential limit $\Sig$ of $\Sig^n$ as in Lemma \ref{lem02}, $(X,L)$ satisfies the SDE \eqref{09} for some $W$, with $x_0=0$.
\end{lemma}

Lemmas \ref{lem02} and \ref{lem03} will be proved in Sections \sect{lem02} and \sect{lem03}, respectively. Before presenting the proofs of Lemmas \ref{lem01}--\ref{lem03}, we show that Theorem \ref{th1} follows.

\noi{\bf Proof of Theorem \ref{th1}.}
By the tightness of $\Sig^n$ stated in Lemma \ref{lem02}, one has relative compactness of their laws. By Lemma \ref{lem03}, the weak limit of $(\hat X^n,\hat I^n)$ along any weakly convergent subsequence of $\Sig^n$ satisfies \eqref{09}. Because Lemma \ref{lem01} states that uniqueness in law holds for solutions to \eqref{09}, it follows that the entire sequence $(\hat X^n,\hat I^n)$ converges weakly to the unique-in-law solution to \eqref{09}.
\qed

\subsection{Tools}

For the proofs of Lemmas \ref{lem01}--\ref{lem03}, we provide some tools regarding martingale problem, $\scrC$-tightness, convergence of martingales.

We first introduce a stopped martingale problem. Let $x \wedge y = \min(x,y)$ for $x,y \in \dd{R}$.
\begin{definition}%[Stopped martingale problem]
\label{dfn:s-mp}
For a relatively open subset $F\subset\R_+$, let $\tau =\tau^F= \inf\{t \ge 0; X(t) \notin F\}$, and define the stopped process $X^{\tau} \equiv \{X^{\tau}(t); t \ge 0\}$ by $X^{\tau}(t) = X(t \wedge \tau)$.
% where $[0,z)$ for $z>0$ is a relative open set on $\R_{+}$.
Let $\scrA$ be the operator defined by \eq{A}. Then, $X$ is said to solve the stopped martingale problem for $(\nu,\scrA, \dd{R}_{+},F)$ if $X(\cdot)=X(\cdot\w\tau)$ a.s., $X(0) \sim \nu$, and
\begin{align}
\label{eq:s-mp}
\begin{split}
  f(X(t))- & \Big(\int_0^{t\w\tau}\scrA f(X(s)) ds + \frac 12 \int_{0}^{t\w\tau} f'(X(s)) dL^{X}_{0}(s)\Big) \\
  & \qquad \mbox{ is an $\dd{F}^{X}$-local martingale for all } f\in {\scrC}^{\infty}_{0}(\R_+,\dd{R}).
  \end{split}
\end{align}
\end{definition}

The martingale problems as well as this stopped ones will be also considered for a different operator or a different domain in the proof of the lemma below. Furthermore, the state space of $X$ may be changed, by which the domain of the corresponding operator is also changed.

The following lemma plays crucial roles in proving Lemma \ref{lem01}. The basic idea of its proof is to construct $X$ by alternatively connecting the independent copies of $X^{\tau_{1}}$ and $X^{\tau_{2}}$. The same idea is used to prove \cite[Section 4.1, Fig.\ 1]{Miya2024b}.

\begin{lemma}
\label{lem:s-mp}
Assume that $X$ is the solution of the martingale problem for $(\nu,\scrA,\dd{R}_{+})$ with $\sr{D}(\scrA)$ for $\scrA$ of \eq{A}. For $z_{1}, z_{2}$ such that $0 < z_{2} < z_{1} < \ell_{1}$, let $F_{1} = [0,z_{1})$ and $F_{2} = (z_{2},\infty)$, and let
\begin{align*}
  \tau_{i} = \inf\{t \ge 0; X(t) \not\in F_{i}\}, \qquad i=1,2.
\end{align*}
Define the stopped process $X^{\tau_{i}} \equiv \{X^{\tau_{i}}(t); t \ge 0\}$ by $X^{\tau_{i}}(t) = X(t \wedge \tau_{i})$ for $i=1,2$. Then, $X^{\tau_{i}}$ is the solution of the stopped martingale problem for $(\nu, \scrA, \dd{R}_{+},F_i)$, which is unique if $X$ is unique, for each $i=1,2$. Furthermore, $X$ is the unique solution of the martingale problem for $(\nu,\scrA,\dd{R}_{+})$ if only if both of $X^{\tau_{1}}$ and $X^{\tau_{2}}$ are the unique solutions of the stopped martingale problems for $(\nu, \scrA, \dd{R}_{+},F_1)$ and $(\nu, \scrA, \dd{R}_{+},F_2)$, respectively.
\end{lemma}

\proof
This lemma can be proved in the same way as \cite[Theorems 4.6.1 and 4.6.2]{eth-kur} for the stopped martingale problems because $F_{1} \cup F_{2} = \dd{R}_{+}$ and $X^{\tau_{i}} \equiv 0$ for $x_0 \not\in F_{i}$ for $i=1,2$.
\qed

\begin{lemma}\label{lem-c}
A sequence of processes $\xi^n$ with sample paths in $\scrD(\R_+,\R)$ is $\scrC$-tight if and only if the following two conditions hold:
\\
For every $T>0$, $\eps>0$, there are $n_0$ and $k$ such that if $n>n_0$,
\[
\PP(\|\xi^n\|^*_T>k)<\eps;
\]
For every $T>0$, $\eps>0$, $\eta>0$, there are $n_0>0$, $\theta>0$ such that if $n>n_0$,
\[
\PP(w_T(\xi^n,\theta)>\eta)<\eps.
\]
\end{lemma}

\proof
See \cite[Proposition VI.3.26, p.\ 351]{jac-shi}.
\qed

\begin{lemma}\label{lem-m1}
For $n\in\N$, let $\check M^n$ be a square-integrable martingale with $\check M^n(0)=0$ a.s. If the sequence $\lan\check M^n\ran$ is $\scrC$-tight in $\scrD(\R_+,\R)$ then the sequence $\check M^n$ is tight in $\scrD(\R_+,\R)$.
\end{lemma}

\proof
See \cite[Theorem VI.4.13, p.\ 358]{jac-shi}.
\qed

\begin{lemma}\label{lem-m2}
For $n\in\N$, let $\check M^n$ be a local martingale, and assume
\[
\sup_n\E(\sup_{s\le t}|\Del\check M^n(s)|)<\iy \qquad \text{for all $t$}.
\]
Then $\check M^n\To M$ implies that $(\check M^n,[\check M^n])\To(M,[M])$.
\end{lemma}

\proof
See \cite[Corollary VI.6.30, p.\ 385]{jac-shi}.
 \qed

\subsection{Proof of Lemma \ref{lem01}.}
\label{sec:lem01}

We first prove that the SDE \eqref{09} has a weak solution $(X,L,W)$. For this, extend $b$ and $\sig$ defined on $\R_{+}$ to $u$  and $v$ defined on $\R$ via
\begin{align}
\label{eq:u-v}
u(x)={\rm sgn}(x)b(|x|), \qquad v(x)={\rm sgn}(x)\sig(|x|), \qquad x \in \R,
\end{align}
where ${\rm sgn}(x)=1$ if $x>0$ and $-1$ if $x\le 0$.
By \cite[Theorem 5.15 p.\ 341]{kar-shr}, the SDE
\begin{align}
\label{eq:QW}
  dQ(t)=u(Q(t))dt+v(Q(t))d\widetilde{W}(t), \qquad Q(0)\sim\widetilde{\nu}.
\end{align}
has a unique weak solution $(Q,\widetilde{W})$, where $\widetilde{W}$ is the standard Brownian motion, and an explosion never occurs in finite time by \cite[Theorem 23.1]{Kall2001} because $1/\sigma^{2}(x)$ is locally integrable. Furthermore, the quadratic variation $[Q](t) = \int_{0}^{t} v^{2}(Q(s)) ds$ is continuous in $t \ge 0$, so $Q$ is a continuous semimartingale by (i) in \sectn{preliminary}.

Define $X(t)=|Q(t)|$ and $W(t) = {\rm sgn}(Q(t)) \widetilde{W}(t)$ for $t \ge 0$, then it follows by Tanaka's formula
\cite[Theorem IV.1.2, p.\ 222]{rev-yor}, $(X,W)$ satisfies the SDE
\begin{align}
\label{eq:XLW}
  dX(t) = b(X(t)) + \sigma(X(t)) dW(t) + dL^{Q}_{0}(t), \qquad x_0 \sim \nu,
\end{align}
where $L_{0}^{Q}$ is the local time of $Q$ at $0$. Now,
\begin{align*}
  L^{Q}_{0}(t) = \lim_{\varepsilon \downarrow 0} \frac 1{2\varepsilon} \int_{0}^{t} 1_{\{- \varepsilon < Q(s) < \varepsilon\}} \sigma_{1}^{2} ds = \lim_{\varepsilon \downarrow 0} \frac 1{2\varepsilon} \int_{0}^{t} 1_{\{0 \le X(s) < \varepsilon\}} d\br{X}(s) = \frac{1}{2}L^{X}_{0}(t),
\end{align*}
because $X(s) = |Q(s)|$ and $\int_{0}^{t} 1(0 \le X(s) < \varepsilon) \br{X}(s) = \sigma_{1}^{2} t$ for $\varepsilon < \ell_{1}$. Hence, $L(t) = L^Q_0(t) = \frac{1}{2}L^X_0(t)$, and therefore the law of $(X,L)$ is determined by that of $X$. Thus, $(X,L,W)$ is the weak solution of the SDE \eqref{09}, equivalently, $(X,W)$ is a weak solution of the SDE \eq{ref-d}. Hence, by \lem{e-mp}, $X$ is also the solution of the martingale problem for $(\nu,\scrA,\dd{R}_{+})$. 

We next show the uniqueness of this solution $X$. To this end, we consider two stopped processes $X^{\tau_{1}}$ and $X^{\tau_{2}}$ defined in \lem{s-mp}. Let $(X_{1},L) \in
\scrC^{\infty}_{[0,\infty)}(\dd{R}_{+}) \times \scrC^\up(\R_+,\dd{R}_{+})$ be the solution of the SDE
\begin{align}
\label{eq:ref-B}
\begin{split}
 & dX_{1}(t) = b_{1} dt + \sigma_{1} dW_{1}(t) + dL(t), \qquad X_{1}(0) \sim \nu,\\
 & \int_{[0,\iy)}X_{1}(t)dL(t)=0. 
\end{split}
\end{align}
where $W_{1}$ is the one dimensional standard Brownian motion. Then, it is well known that $X_{1}$ is the unique strong solution of the SDE \eq{ref-B}, called the reflecting Brownian motion on $[0,\infty)$. Furthermore, $X_{1}$ is a reflecting process similar to the solution $X$ of the SDE \eqref{09}, and therefore $L = \frac{1}{2}L^{X_{1}}_{0}$ similarly by (iii) in \sectn{preliminary}. 

Define operator $\scrA_{1}$ with domain ${\scrC}^{\infty}_{0}(\R_+,\dd{R})$ as
\begin{align}
\label{eq:Af}
  \scrA_{1} f(x)=b_{1} f'(x)+\frac{1}{2} \sigma_{1}^{2} f''(x), \qquad f \in {\scrC}^{\infty}_{0}(\R_+,\dd{R}).
\end{align}
 Since the SDE \eq{ref-B} has the unique solution and $L = \frac{1}{2}L^{X_{1}}_{0}$, $X_{1}$ is the unique solution of the martingale problem $(\nu,\scrA_{1},\dd{R}_{+})$ by \lem{e-mp} for $\scrA = \scrA_{1}$. Hence, by \lem{s-mp}, the stopped process $X_{1}^{\tau_{1}}$ is also the unique solution of the stopped martingale problems for $(\nu,\scrA_{1},\dd{R}_{+},F_1)$. Furthermore, $X_{1}^{\tau_{1}}$ is identical with $X^{\tau_{1}}$ in law. Since we already see that $X$ is the solution of the martingale problem for $(\nu,\scrA,\dd{R}_{+})$, $X^{\tau_{1}}$ is the unique solution of the stopped martingale problem for $(\nu,\scrA,\dd{R}_{+},F_1)$ by \lem{s-mp}.

We next consider the process $Q$, and recall that it is the unique weak solution of the SDE \eq{QW}. Then, $Q$ is the unique solution of the martingale problem for $(\nu,\scrG,\R)$ by \cite[Theorem 21.7]{Kall2001}, where $\scrG$ is given by \eq{G} with $u, v$ of \eq{u-v}. Hence, $Q^{\tau_{2}}$ is also the unique solution of the stopped martingale problem for $(\nu,\scrA,\dd{R},F_2)$ by \lem{s-mp}. Furthermore, $Q^{\tau_{2}}$ is identical with $X^{\tau_{2}}$ in law because $\scrA f(x)$ for $x \in F_{2}$ is unchanged and both of $Q^{\tau_{2}}$ and $X^{\tau_{2}}$ are positive valued, so $X^{\tau_{2}}$ is the unique solution of the stopped martingale problem $(\scrA,\nu,\dd{R}_{+},F_2)$.

Thus, $X^{\tau_{1}}$ and $X^{\tau_{2}}$ are the unique solutions of the stopped martingale problems for $(\nu,\scrA,\dd{R}_{+},F_1)$ and $(\nu,\scrA,\dd{R}_{+},F_2)$, respectively. Hence, by \lem{s-mp}, $X$ is the unique solution of the martingale problem $(\nu,\scrA,\dd{R}_{+})$, and therefore $X$ is also the unique solution of the SDE \eq{ref-d} by \lem{e-mp}.
\qed

\subsection{Proof of Lemma \ref{lem02}.}
\label{sec:lem02}

The proof proceeds in several steps.

Step 1.
The sequence $(\bar A^n,\bar D^n)$ is $\scrC$-tight,
and, for each $t$, the sequences $\bar A^n(t)$
and $\bar D^n(t)$ are uniformly integrable.

To prove the claim regarding $\scrC$-tightness, recall that $A$ is a renewal process
for which the first moment of inter-event distribution has mean $1$.
Therefore $n^{-1}A(n\cdot)\to\iota$ in probability. Let $\bar A^n(t)= n^{-1} A^n(t)$, then
\begin{equation}\label{12}
\bar A^n(t)=n^{-1}A\Big(n\sum_{i=0}^K\bar\la^n_iH^n_i(t)\Big),
\end{equation}
and noting that $\bar\la^n_i$ are bounded and $H^n_i$
are 1-Lipschitz, gives $\scrC$-tightness of $\bar A^n$ in view of Lemma \ref{lem-c}.
A similar conclusion holds for $\bar D^n$.

As for the claim of uniform integrability, in view of \eqref{12}, the boundedness of $\bar\la^n_i$ and the bound $H^n_i(t)\le t$,
one has for fixed $t$, $\bar A^n(t)\le n^{-1}A(cn)$ for some deterministic $c=c(t)$.
Uniform integrability of $n^{-1}A(cn)$ is well known
(see \cite[(5.6), Chapter 2, page 58]{gut09}),
proving the claim for $\bar A^n(t)$, and similarly for $\bar D^n(t)$.

Step 2.
Let $Z(j)$, $j=1,2,\ldots$ be a collection of i.i.d.\ random variables such that $\E[Z(1)^2]<\iy$. Then
\begin{equation}\label{101}
n^{-1/2}\max_{j\le n}Z(j)\to0 \text{ in probability, as } n\to\iy.
\end{equation}
This is a standard fact, but for completeness we provide
a proof. Fix $\eps>0$. Then, denoting $Z_\eps=\eps^{-2}Z(1)^2$,
\begin{align*}
\PP(n^{-1/2}\max_{j\le n}Z(j)>\eps)
&\le n\PP(Z(1)>\eps n^{1/2})
\le n\PP(Z(1)^2>\eps^2n)
\\
&=2\frac{n}{2}\PP(Z_\eps>n)
\le 2\sum_{[n/2]\le k\le n}\PP(Z_\eps>k)
\le 2\sum_{k\ge [n/2]}\PP(Z_\eps>k).
\end{align*}
The last expression converges to $0$ as $n\to\iy$,
because, by assumption, $\E[Z_\eps]<\iy$. This proves
\eqref{101}.

Step 3. The processes $e^n$, $e^n_A$ and $e^n_{S}$ converge
to $0$ in probability in $\scrD(\R_+,\R)$.

We prove the claim regarding $e^n_A$
(a similar proof holds for $e^n_S$).
By Lemma \ref{lem1},
\[
e^n_A(t)=n^{-1}\sum_{j=1}^{A^n(t)}\beta(j),
\qquad
\beta(j)=\zeta_A(j)^2-\sig^2_A.
\]
Denoting $\gamma(N)=\sum_{j=1}^N\beta(j)$, $N\in\Z_+$,
we have
$e^n_A(t)=n^{-1}\gamma(A^n(t))$.
Fix $t$. Given $\eps$, using Step 1, we can find
$c$ such that $\PP(\bar A^n(t)>c)<\eps$.
On the event $\bar A^n(t)\le c$,
\[
\|n^{-1}\gamma(A^n(\cdot))\|^*_t
\le n^{-1}\|\gamma\|^*_{[cn]}.
\]
Now, $\beta(j)$ are i.i.d.\ with mean 0. By the functional LLN,
$N^{-1}\|\gamma\|^*_N\to0$ in probability as $N\to\iy$. Thus for any $\del>0$
\[
\PP(\bar A^n(t)\le c,\|e^n_A\|^*_t>\del)
\le\PP(n^{-1}\|\gamma\|^*_{[cn]}>\del)\to0,
\]
as $n\to\iy$. Hence $\limsup_n\PP(\|e^n_A\|^*_t>\del)\le \eps$.
Since $\eps>0$ is arbitrary, this shows that $e^n_A\to0$ in probability.

We now show that $\|e^{n}\|^{*}_{t} \to 0$ in probability as $n \to \infty$ for each fixed $t \ge 0$.
It suffices to show
\begin{align}\label{a01}
  \sup_{0 \le s \le t} n^{-1/2}R_A\left(\sum_{i=0}^K\la^nH^n_i(s)\right) \to 0 \; \mbox{in probability as $n \to \infty$},
\end{align}
and a similar statement for the term involving
$R_S$ in \eqref{13}.
To this end, using \eqref{05}, we can write
\[
\sup_{0\le s\le T} R_A(s) = \max_{0\le j\le A(T)}Z_A(j).
\]
Fix $c>0$, $\eps>0$ and $\del>0$.
Let $n_0$ and $c_1<\iy$ be such that
$\PP(A(cn)>c_1n)<\del$ for
all $n>n_0$ (such constants exist again by the fact that
$n^{-1}A(n\cdot)\to\iota$ in probability). Then
\begin{align*}
\PP(n^{-1/2}\|R_A\|^*_{cn}>\eps)
&\le
\PP(A(cn)>c_1n)+\PP(n^{-1/2}\max_{0\le j\le c_1n}Z_A(j)>\eps).
\end{align*}
The first term on the right is bounded by $\del$ for $n>n_0$,
and the second term converges to $0$ as $n\to\iy$,
by Step 2. Sending $\del\to0$, it follows that
$\lim_{n\to\iy}\PP(n^{-1/2}\|R_A\|^*_{cn}>\eps)=0$. In other words, $n^{-1/2}\|R_A\|^*_{cn}\to0$ in probability.
for any constant $c$. Because this constant $c$ is arbitrary,
$\sum_iH^n_i(t)= t$, and $\bar\la^n_i$ are bounded, this shows that \eqref{a01} holds.
A similar argument holds for the term involving $R_S$,
and it follows that $\|e^{n}\|^{*}_{t} \to 0$ for each fixed $t \ge 0$.

Step 4.
The sequences $\lan\hat M^n\ran$ and $[\hat M^n]$ are $\scrC$-tight
in $\scrD(\R_+,\R)$.

In view of Lemma \ref{lem1}, the $\scrC$-tightness of
$\lan\hat M^n_A\ran$ and $\lan\hat M^n_S\ran$ follows from that of $\bar A^n$
and $\bar D^n$, proved in Step 1. In view of \eqref{42},
\[
\lan\hat M^n\ran=\lan\hat M^n_A\ran -2\lan\hat M^n_A,\hat M^n_S\ran
+\lan\hat M^n_S\ran
=\lan\hat M^n_A\ran+\lan\hat M^n_S\ran,
\]
and the $\scrC$-tightness of $\lan\hat M^n\ran$ follows.

By Step 3 and the $\scrC$-tightness of $\lan\hat M^n_A\ran$ and $\lan\hat M^n_S\ran$, $[\hat M^n_A]$ and $[\hat M^n_S]$ are also $\scrC$-tight. To prove
$\scrC$-tightness of $[\hat M^n]$ it thus suffices to prove $\scrC$-tightness
of $m^n=[\hat M^n_A,\hat M^n_S]$. We argue as follows. Given $T$ and $\del$,
\begin{align*}
w_T(m^n,\del)&=\sup\Big\{\Big|
\sum_{t_1<s\le t_2}\Del\hat M^n_A(s)\Del\hat M^n_S(s)
\Big|:0\le t_1<t_2\le T, \, t_2-t_1\le \del\Big\}
\\
&\le\sup\Big\{
\sum_{t_1<s\le t_2} \frac 12 (\Del\hat M^n_A(s)^2+\Del\hat M^n_S(s)^2)
:0\le t_1<t_2\le T, \, t_2-t_1\le \del\Big\}
\\
&\le \frac 12 (w_T([\hat M^n_A],\del)+w_T([\hat M^n_S],\del)).
\end{align*}
Similarly, $\|m^n\|^*_T\le ([\hat M^n_A](T)+[\hat M^n_S](T))/2$.
As a result, $\scrC$-tightness of $m^n$ follows from that of
$[\hat M^n_A]$ and $[\hat M^n_S]$, by virtue of Lemma \ref{lem-c}.

Step 5. The sequence $\hat M^n$ is tight in $\scrD(\R_+,\R)$.
This is immediate from Step 4, using Lemma~\ref{lem-m1}.
 
Step 6. The sequence $(\hat M^n,[\hat M^n])$
is $\scrC$-tight, and
if, along a subsequence, $\hat M^n\To M$, then
 $(\hat M^n,[\hat M^n])\To(M,[M])$ along this subsequence.

To prove this statement, we first verify the condition
\begin{equation}\label{14}
\sup_n\E(\sup_{s\le t}|\Del\hat M^n(s)|)<\iy \qquad \text{for all $t$}.
\end{equation}
To this end, write
\[
\E(\sup_{s\le t}|\Del\hat M^n_A(s)|)
\le 1+\int_1^\iy\PP(n^{-1/2}\sup_{j\le A^n(t)}|\zeta_A(j)|>r)dr
\le 1+\int_1^\iy\frac{\E(\sup_{j\le A^n(t)}\zeta_A(j)^2)}{nr^2}dr.
\]

For $t\ge1$, we have $A^n(t)\le A(cnt)$ for suitable constant $c$,
and so
\[
\E\Big[\sup_{j\le A^n(t)}\zeta_A(j)^2\Big]
\le \E\Big[\sup_{j\le A(cnt)}\zeta_A(j)^2\Big]
\le\E\Big[\sum_{j=1}^{A(cnt)}\zeta_A(j)^2\Big]
=\E(A(cnt))\E(\zeta_A(1)^2).
\]
Since $\E(A(cnt))\le c_1nt$
for some constant $c_1$, we obtain that $\E(\sup_{s\le t}|\Del\hat M^n_A(s)|)$
is bounded uniformly in $n$. A similar statement holds for $\hat M^n_S$,
and, since $|\Del\hat M^n|\le|\Del\hat M^n_A|+|\Del\hat M^n_S|$,
also for $\hat M^n$. This verifies \eqref{14}.

Consider a subsequence along which $\hat M^n$ converges in distribution
in $\scrD(\R_+,\R)$ and denote its limit by $M$.
We can now apply Lemma \ref{lem-m2}
by which, under condition \eqref{14},
$(\hat M^n,[\hat M^n])\To(M,[M])$ holds along the same sequence. Next,
in view of Step 4, $[M]$ is continuous. This implies that
$M$ is also continuous by item (i) of \sectn{preliminary}.
Therefore $(\hat M^n,[\hat M^n])$ is $\scrC$-tight.

Step 7. $\Sig^n$ is $\scrC$-tight.

To prove this statement, it remains to show the $\scrC$-tightness of $\hat X^n$, $\hat I^n$, $\hat Y^n$ and $H^n$.
To this end, note that $\hat X^{n} = \hat Y^{n} + \hat I^{n} \ge 0$ and, by definition,
\[
\int_{[0,\iy)}\hat X^n(t)d\hat I^n(t)=0.
\]
As a result, we can rewrite \eqref{10} as
\begin{equation}
\label{17}
(\hat X^n,\hat I^n)=\Gamma(\hat Y^n).
\end{equation}
Now, $H^n_i$ are uniformly Lipschitz and are therefore $\scrC$-tight.
Also, $b^n_i$ are bounded,
therefore it follows from Steps 1 and 6 that $\hat Y^n$ are $\scrC$-tight.
The Lipschitz continuity of the map $\Gamma:
\scrD(\R_+,\R)\to \scrD(\R_+,\R)^2$ with respect to the uniform
topology on compacts, implies that the sequence \eqref{16}
$\scrC$-tight, and that the relation \eqref{17} is preserved under the limit.
Moreover, the tightness of $\hat I^n$ implies that $H_0=0$.

Step 8.
Fix a convergent subsequence of \eqref{16} and denote its weak limit by \eqref{18}. Then $M$ is a martingale and relations \eqref{20-} and \eqref{20} hold.

Because $M$ is the weak limit of the martingales $\hat M^n$,
to show that $M$ is a martingale, it suffices to show that
for fixed $t$, $\hat M^n(t)$ are uniformly integrable.
This is true by
$\sup_n\E(\hat M^n(t)^2)<\iy$, which follows from Lemma \ref{lem1}
and $\sup_nn^{-1}\E A^n(t)<\iy$. Hence $M$ is a martingale.

As for \eqref{20-}, we have already shown its first two parts. Its last part is immediate from \eqref{10}.

Next, we show \eqref{20}. The fact that the limit $\bar A$ of $\bar A^n$ is given by the expression in \eqref{20} follows from \eqref{12} and the convergence $n^{-1}A(n\cdot)\to\iota$. The form \eqref{20} of $\bar D$ follows similarly.

It remains to show that $[M]$ is given by the expression in \eqref{20}.
By Lemma \ref{lem1} and the convergence of $(\bar A^n,\bar D^n)$, we have
$(\lan \hat M^n_A\ran,\lan \hat M^n_S\ran)\To(\sig^2_A\bar A,\sig^2_S\bar D)$.
Because $\lan\hat M^n_A,\hat M^n_S\ran=0$, it follows that
$\lan\hat M^n\ran\To\sig^2_A\bar A+\sig^2_S\bar D$, which is $\sum_{i=1}^K\sig^2_iH_i$, in
view of $\sig_i^2=\la_i\sig_A^2+\mu_i\sig_S^2$, $1\le i\le K$.

By Step 6, $(\hat M^n,[\hat M^n])\To(M,[M])$ along the
convergent subsequence. Suppose we also show that the limit along this subsequence satisfies
\begin{equation}\label{102}
\lan M\ran=\tilde H:=\sum_{i=1}^K\sig_i^2H_i.
\end{equation}
Since $M$ is a continuous martingale, $[M]=\lan M\ran$,
and thus the expression for $[M]$ in \eqref{20} would be proved.

Thus, it remains to show \eqref{102}.
To this end, note that we have shown that
\[
\La^n:=(\hat M^{n})^{2} - \lan\hat M^n\ran \Rightarrow \La:= M^{2} - \tilde H.
\]
Let us show that $\La^n(t)$ are uniformly integrable for each $t$.
Recall from Step 1
that $\bar A^n(t)$ and $\bar D^n(t)$ are uniformly integrable
for each fixed $t$. In view of (\ref{42}), this implies that $\lan\hat M^n\ran(t)$ is uniformly integrable. By the definition \eqref{08},
we can apply
\cite[Theorem 6.2, Chapter 1, page 32]{gut09} to $\hat M^{n}_{A}$ and $\hat M^{n}_{S}$, which gives
the uniform integrability of $(\hat M^{n}_{A})^{2}(t)$ and $(\hat M^{n}_{S})^{2}(t)$. This proves the uniform integrability of $(\hat M^{n})^{2}(t)$ because $(\hat M^{n})^{2}(t) \le 2[(\hat M^{n}_{A})^{2}(t) + (\hat M^{n}_{S})^{2}(t)]$. Hence, $\La^n(t)$ is uniformly integrable.
Next,
let $\tilde{\sr{F}}_t=\sr{F}^{\tilde H}_t\vee\sr{F}^M_t$. We will now show that
$\La$ is an $\{\tilde{\sr{F}}_t\}$-martingale. Since this process
is adapted and integrable, it remain to show that for $0\le s\le t$,
$\E[\La(t)-\La(s)|\tilde{\sr{F}}_s]=0$ a.s.
For this, it suffices that for any bounded continuous
$h:\scrD(\R_+,\R)\times \scrC(\R_+,\R)\to\R$,
\begin{equation}\label{103}
\E[h(M(\cdot\w s),\tilde H(\cdot\w s))(\La(t)-\La(s))]=0.
\end{equation}
Since $\La^n$ is an $\{\sr{F}^{n}_{t}\}$-martingale for each $n \ge 1$,
we have
\[
\E[h(\hat M^n(\cdot\w s),\lan \hat M^n\ran(\cdot\w s))(\La^n(t)-\La^n(s))]=0.
\]
The convergence $(\hat M^n,\lan\hat M^n\ran,\La^n)\To(M,\tilde H,\La)$,
the boundedness of $h$ and uniform integrability of $\La^n(s)$
and $\La^n(t)$ imply \eqref{103}.
This shows that $\La$ is an $\{\tilde{\sr{F}}_t\}$-martingale.
Since $\tilde H$ is adapted to this filtration and continuous,
it is predictable on it.
Hence by the uniqueness of the Doob-Meyer decomposition for the submartingale $M^{2}$, we conclude that
$\lan M \ran = \sum_{i=1}^K\sig^2_iH_i$.
This proves \eqref{102} and completes the proof of the lemma.
\qed

\subsection{Proof of Lemma \ref{lem03}.}
\label{sec:lem03}

Fix a convergent subsequence of $\Sig^n$ (as in \eqref{16}), and denote its limit by $\Sig$ (as in \eqref{18}). The convergence of this tuple occurs in a space $\scrD(\R_+,\R)^{10+K}$, a Polish space. Thus one can
invoke Skorohod's representation theorem \cite[Theorem 25.6]{bill}, and assume without loss of generality that this convergence is a.s. We show that this a.s.\ limit establishes \eqref{09}, which proves Lemma \ref{lem03}.

Denote $\sr{F}_t=\sr{F}^X_t\vee\sr{F}^M_t$.
We show that $M$ is an $\{\sr{F}_t\}$-martingale,
by an argument similar to the one in Step 8 of the proof of Lemma \ref{lem02}.
Clearly, it is adapted and integrable, so it remains to show that for $0\le s\le t$,
$\E[M(t)-M(s)|\sr{F}_s]=0$ a.s. For this, it suffices
to show that for any bounded continuous $h:\scrD(\R_+,\R)^2\to\R$,
\begin{equation}\label{22}
\E[h(X(\cdot\w s),M(\cdot\w s)) (M(t)-M(s))]=0.
\end{equation}
Since $\hat M^n$ is an $\{\sr{F}^n_t\}$-martingale, we have
\[
\E[h( \hat X^n(\cdot\w s) ,\hat M^n(\cdot\w s) ) (\hat M^n(t)-\hat M^n(s))]=0.
\]
Moreover, $(\hat X^n,\hat M^n)\to(X,M)$ a.s.\ in $\scrD(\R_+,\R)^2$,
and so \eqref{22} will follow provided that
$h(\hat X^n(\cdot\w s) ,\hat M^n(\cdot\w s) ) (\hat M^n(t)-\hat M^n(s))$ are uniformly integrable. This last condition
holds by the boundedness of $h$ and the fact
that $\sup_n\E(\hat M^n(t)^2)<\iy$, which follows from Lemma \ref{lem1}
and $\sup_nn^{-1}\E A^n(t)<\iy$.

Next we show that
\begin{equation}\label{15}
\int_0^\iy1_{\{X_t=\ell\}}dt=0,\qquad \ell=\ell_0,\ell_1,\ldots,\ell_{K-1},
\qquad {\it a.s.}
\end{equation}
To this end, note first that, a.s., $d\lan M\ran(t)=d[M](t)\ge c_1dt$ as measures, where
$c_1=\min_{1\le i\le K}\sig^2_i>0$.
This follows from the fact that $\sum_{i=1}^KH_i(t)=t$ (recall that $H_0=0$),
and the expression for $[M]$ in \eqref{20}.

Fix $0\le i\le K-1$ and set $\ell=\ell_i$.
Then by the occupation time formula \cite[Corollary VI.1.6]{rev-yor},
%with $L^{X}_{a}(t)$ the local time of $X$ at $a\in[0,\iy)$ by time $t$,
\[
\int_0^t1_{\{X_s=\ell\}}d\lan X\ran(s)
=\int_0^\iy L^{X}_{a}(t) 1_{\{a=\ell\}} da=0.
\]
But $d\lan X\ran(t)=d\lan M\ran(t)\ge c_1dt$, and \eqref{15} follows.

Finally, we show that for any subsequential limit \eqref{18} of \eqref{16},
the law of $X$ is a solution to the martingale problem.
% In view of Lemma \ref{lem2}, this proves Claim 3.

%The convergence $\hat X^n\to X$ along with Step 9
%imply that $H_i$, $1\le i\le K$, which are the limits of
%$H^n_i=\int_0^\cdot 1_{\{\hat X^n(s)\in\bS_i\}}ds$, are given by
%$H_i=\int_0^\cdot 1_{\{X(s)\in\bS_i\}}ds$.

The convergence $\hat X^n\to X$ and the continuity of $X$ imply that for fixed $t>0$, $\hat X^n\to X$ uniformly in $[0,t]$ a.s. As a result, $\int_0^t g(\hat X^n(s))ds\to\int_0^tg(X(s))ds$ a.s.\ for every $g$ that is continuous and bounded. For any deterministic sequence $c_n\to0$, we similarly have
\[
\int_0^t g(\hat X^n(s)+c_n)ds\to\int_0^tg(X(s))ds
\]
a.s.\ for all such $g$. Thus $\tilde m_n\to \tilde m$ a.s., where $\tilde m_n$ and $\tilde m$ are the probability measures on $\R$ defined by
\[
\tilde m_n([a,b])=t^{-1}\int_0^t1_{\{\hat X^n(s)+c_n\in[a,b]\}}ds,
\qquad
\tilde m([a,b])=t^{-1}\int_0^t1_{\{X(s)\in[a,b]\}}ds,
\qquad -\iy<a<b<\iy.
\]
Hence Portmanteau's theorem gives that, for $C>0$, if $\tilde m(\{C\})=0$ a.s.\ then a.s.,
\[
\tilde m_n((-\iy,C])\to \tilde m((-\iy,C]).
\]
To use this fact, fix $1\le i\le K-1$ and let $C=\ell_i$ and $c_n=-n^{-1/2}\ell^n_i+\ell_i$. By \eqref{eq:ell-n}, $c_n\to0$. By \eqref{15}, $\tilde m(\{C\})=0$ a.s. This shows that, a.s.,
\[
\int_0^t1_{\{X^n(s)\le\ell^n_i\}}ds
=\int_0^t1_{\{\hat X^n(s)+c_n\le\ell_i\}}ds\to\int_0^t1_{\{X(s)\le\ell_i\}}ds.
\]
From this, we obtain that $H_i$, $1\le i\le K$, which are the limits of
$H^n_i=\int_0^\cdot 1_{\{X^n(s)\in\bS^n_i\}}ds$, are given by
$H_i=\int_0^\cdot 1_{\{X(s)\in\bS_i\}}ds$.
As a result,
\[
\sum_{i=1}^Kb_iH_i(t)=\int_0^tb(X(s))ds,
\]
where the expression \eqref{11} for $b$ is used. Hence, let $Y(t) = \int_0^tb(X(s))ds+M(t)$, then, from \eqref{20-}, we have
\begin{align}
\label{eq:Xn-limit}
  X(t)=\int_0^tb(X(s))ds+M(t)+L(t), \qquad (X,L) = \Gamma(Y).
\end{align}

As for $[M]$, the expression in \eqref{20} and the fact that $H_0=0$ shows that
\[
[M](t)=\sum_{i=1}^K\sig^2_iH_i=\sum_{i=0}^K\sig^2_iH_i
=\int_0^t\sig^2(X(s))ds,
\]
where the expression in \eqref{11} for $\sig$ is used.
As explained in the proof of \lem{e-mp} in Appendix, this implies that there exists a standard Brownian motion $W$ which is independent of $X$ such that
\begin{align*}
  M(t) = \int_{0}^{t} \sigma(X(s)) dW(s), \qquad t \ge 0.
\end{align*}
Hence, \eq{Xn-limit} implies \eqref{09}.
\qed

%\newpage

\appendix
\section*{Appendix}
\beginsec
\setcounter{section}{1}

\subsection{Proof of \lem{e-mp}}
\label{app:e-mp}

Assume that a continuous semimartingale $X$ is the solution of the SDE \eq{ref-d}. Then, applying Ito's formula
to \eq{ref-d} for test function $f\in {\scrC}^{\infty}_{0}(\R_+,\dd{R})$,
\begin{align}
\label{eq:fXLW}
\begin{split}
f(X(t)) & = f(x_0)+\int_0^t \left(\scrA f(X(s))ds+\frac 12 f'(X(s)) dL^{X}_{0}(s)\right)\\
& \qquad  + \int_0^t f'(X(s)) \sig(X(s))dW(s).
\end{split}
\end{align}
This equation shows that $X$ solves the martingale problem $(\nu,\scrA,\dd{R}_{+})$ because
\begin{align*}
  \int_0^t f'(X(s)) \sig(X(s))dW(s) \mbox{ is a martingale}. 
\end{align*}
Conversely, assume that $X$ solves the martingale problem for $(\nu,\scrA,\R_{+})$. For $i=1,2$, denote the local martingale in \eq{e-mp} for $f(x) \in {\scrC}^{\infty}_{0}(\R_+,\dd{R})$ such that $f(x) = x^{i}$ for $x \in [0,n)$ for each fixed integer $n \ge 1$ by $M_{i} \equiv \{M_{i}(t); t \ge 0\}$, and let $\chi(t) = \int_0^t 2 b(X(s)) ds + L^{X}_{0}(t)$. Since $\chi$ is continuous and has bounded variations, we have, through localization arguments,
\begin{align}
\label{eq:M1}
 M_{1}(t) & = X(t) - x_0 - \frac 12 \chi(t), \qquad t \ge 0,\\
\label{eq:M2}
 M_{2}(t) & = X^{2}(t) - X^{2}(0) - \int_0^t (X(s) d\chi(s) + \sigma^{2}(X(s)) ds), \qquad t \ge 0,
\end{align}
and both of $M_{1}$ and $M_{2}$ are continuous. On the other hand, from Ito integration formula to $X$ for test function $f(x)=x^{2}$ and \eq{M1},
\begin{align*}
  X^{2}(t) & = 2 \int_{0}^{t} X(s) dX(s) + [X](t)\\
  & = 2 \int_{0}^{t} X(s) dM_{1}(s) + \int_{0}^{t} X(s) \chi(ds) + [X](t), \qquad t \ge 0,
\end{align*}
where $[X]$ is the quadratic variation of $X$. Substituting this $X^{2}(t)$ into \eq{M2} and writing $dX(s)$ as $dM_{1}(s) + 2^{-1} d \chi(s)$ by \eq{M1}, we have
\begin{align*}
  M_{2}(t) - 2 \int_{0}^{t} X(s) dM_{1}(s) = [X](t) - \int_{0}^{t} \sigma^{2}(X(s)) ds.
\end{align*}
The left-hand side of this equation is a continuous martingale, while its right-hand side is a process with bounded variations. Since a continuous martingale cannot have a bounded variation component as is well known, both sides must vanish, and therefore, using \eq{M1}, we have
\begin{align}
\label{eq:M0}
  [M_{1}](t) = [X](t) = \int_0^t \sig^{2}(X(s))ds, \qquad t \ge 0.
\end{align}
Thus, by the well known characterization of Brownian motion by a local martingale (e.g., see \cite[Theorem 18.12]{Kall2001}), there exists a standard Brownian motion $\widehat{W}$ which is independent of $X$ such that
\begin{align*}
  M_{1}(t) = \int_{0}^{t} \sigma(X(s)) d\widehat{W}(s), \qquad t \ge 0,
\end{align*}
 This together with \eq{M1} implies the SDE \eq{ref-d}. Thus, if $X$ is the solution of the martingale problem for $(\nu,\scrA,\R_{+})$, then $(X,\widehat{W})$ is the weak solution of the SDE \eq{ref-d}. This together with the first paragraph completes the proof of this lemma.
 \qed

\medskip

{\bf Acknowledgment.}
The authors thank an associate editor and two referees for their constructive suggestions, which have significantly improved the exposition.
The first author was supported by ISF grants 1035/20 and 3240/25.

\footnotesize

%\bibliographystyle{is-abbrv}
%\bibliography{../../../../texmf/bib/dai20230708M3}
%\end{document}

%%\bibliographystyle{plain}
%%\bibliographystyle{annotate}
%%\bibliographystyle{apalike}
\bibliographystyle{is-abbrv}

\bibliography{main-rev}

%\begin{thebibliography}{90}
%
%
%\end{thebibliography}

\end{document}